\documentclass[12pt,letterpaper]{article}

\usepackage[margin=1in]{geometry} 

\usepackage{amsmath}
\usepackage{amsthm}
\usepackage{amssymb}

\usepackage[utf8]{inputenc}
\usepackage{hyperref}
\hypersetup{colorlinks,allcolors=black}

\theoremstyle{plain}

\theoremstyle{definition}

\usepackage{graphicx, color}

\usepackage{algorithm, algpseudocode} 
\usepackage{mathrsfs} 

\title{ Disentangling global equity market instability:\\ a network analysis}
\author{Supanat Kamtue \and Pongsak Luangaram \and Sirawit Woramongkhon\thanks{Electronic address: \texttt{sirawitw@bot.or.th}: Corresponding author}}

\usepackage[margin=1in]{geometry}
\usepackage{amsmath, amssymb, amsthm}
\usepackage{graphicx}

\usepackage{enumerate}
\usepackage{caption}
\usepackage{subcaption}

\usepackage{xcolor}

\usepackage{bm}
\usepackage{comment}

\usepackage[
backend=biber,
style=alphabetic,
sorting=ynt
]{biblatex}
\begin{document}

	\maketitle
    \begin{abstract}
During a financial crisis, the capital markets network frequently exhibits a high correlation between returns.
We developed a network analysis framework based on daily returns from 42 countries to determine systemic stability.
Our network is built using the conditional probability of co-movement of returns, and it identifies nodes, network complexity, and edge as potential sources of fragility.
We also introduce the concept of measuring flows from one return to another.
Then, we use 120-day rolling data to capture the financial system's behavior and create a financial stability indicator.
We discover that the contributions of nodes and network complexity to changes in system stability frequently cancel each other out. Edge change may be a determinant of systemic stability.
Furthermore, the total flows in the network are highly correlated with the volatility. 
It main advantage is the tractability and potential sources of volatility can be determined.  
    \end{abstract}
    
\section{Introduction}

The complex financial system is continuously evolving, as does the network architecture culminating in systemic vulnerability.
A single collapse in a fragile structure can precipitate subsequent ruptures throughout the network, greatly amplifying the magnitude of the impact.
Therefore, understanding systemic transition into instability is critical for assessing market behavior, optimizing portfolios and formulating financial policy.
Quantifying an early warning signal or indication of the stability in the financial network can navigate future expectations of market performance and  call for necessary adjustments to prepare for and possibly avoid a full-fledged collapse.
Scheffer et al \cite{stability-intro} reviews the framework for determining indicators for such critical transitions;
because financial system is  highly volatile and stochastic nature of the financial system, the variance of asset returns can be indicative of the network's flickering. 
It can then be used to identity the probability of entering unstable states. 
The goal of this article is to introduce quantitative method, based on the concept of network flows, to construct an indicator that could both reflect network fragility and provide proper criterion for instability period.
We then create this indicator from capital market returns across 42 countries. 

The financial system stability is defined in terms of the network's ability to absorb shock contagion caused by one or more agents.
In an unstable financial system, a shock emanating from one country can precipitate the collapse of multiple countries.
The existing literature attempts to quantify systemic risk and develop an early warning indicator of financial market fragility and instability through the analysis of Ricci curvature of the financial network.
The development of risk measures are both on specific financial assets or institutions (nodes) and, their interrelationships (edges).
In \cite{ricci}, \cite{ricci2} and \cite{ricci3}, the authors provide empirical evidence suggesting that, by analyzing a financial system as a weighted network, the average of Ricci-type curvature may serve as an indicator for systemic fragility.
Despite such a development, to our knowledge, one disadvantage of their approach is that there is no known criterion for separating stability and instability periods.

Our goal is to introduce quantitative method, based on the concept of network flows, to construct an indicator that could both reflect network fragility and to provide proper criterion for instability period.
Our flow quantifier is heavily inspired by inbound and outbound-edge Forman Ricci curvature which appears in the formula (3)-(5) in \cite{sreejith2017forman}. 
The analysis starts with the construction of directed network of flow between two capital markets based on conditional probability of concurrent abruption. 
The incoming and outgoing flow of a node can then roughly represents the probabilistic impact of the respective market performance.
Afterwards, we design the flow movement process therein and determine the largest eigenvalue of the underlying transitional matrix.

Our major contribution is the development and derivation of the fragility indicator based on flow process and spectrum that could serve as criterion distinguishing period of instability. 
In fact, the network is unstable when the associated spectrum is higher than $1$, and stable otherwise. 
We demonstrate that our analysis on the flow movement could provide not only the systemic fragility but also general behavior and potential internal instability of network's components as well.
We analyze the the ratio of total incoming and outgoing flows and its distribution and discover their association with instability in the respective country and relationship. 
Additionally, the movement of total systemic flows greatly share similarity with the market volatility. 
The analysis on the major components of total flows in each time period can roughly represent the sources of market volatility.  

This paper is organized as follows.
Section 2 describes the data used in our analysis.
Then, it introduces probabilistic flow concepts and shock transmission process based on conditional probability. 
Afterwards, it outlines stability indicator calculation and the underlying contribution to change. 
Section 3 demonstrates our results in the application of network flows in capturing market behavior.
It also discuss to performance of our financial fragility indicators in showing unstable periods in the global financial network.
Then, it analyze the main sources if financial instability from the network's features.
Further more, the potential connection to market volatility is discussed.
Section 4 concludes.

\section{Methodology}

\subsection{Data and preparation}

The data set consists of 5-business-day MSCI indexes from January 1995 to July 2021 for $n=42$ countries.
Each index is calculated from equities in large and medium-sized companies from a variety of industries.
The underlying construction of indexes implies that the percentage growth of the MSCI index at period $t$ or $r_t$ can reflect performance of the capital market in the respective country; \begin{equation}
\label{def:return}
	r_t = \frac{i_t}{i_{t-\triangle t}}-1, 
\end{equation}
where $\triangle t$ denotes time-shift and $i_t$ is the index at time $t> \triangle t.$ 
Our analysis sets $\triangle t = 1$ to analyze daily return.
The log difference method is an alternative to our simple return formula (\ref{def:return}).
Both approaches deliver comparable and nearly identical value when the financial market does not experience tailed events and the index grows slowly.
However, the log difference return may underestimate the stock market performance when the index grows abruptly while over-exaggerating during collapses or crisis.
Thus, the simple return approach may provide an accurate picture of stock market performance during tail events.

The market movement can be captured via daily return difference $d_t = r_t-r_{t-1}$.
The positive (negative) value of $d t$ indicates that the stock market performance is better (worse) than the previous day's.
The extreme tailed events such as financial collapse or pandemic outbreak will likely result in negative value of $d_t$. 
The return difference has the advantage that it can capture the change in stock market performance at the start of the tailed events.
For example, the stock market may have performed well the day before the crisis, with a positive return.
When the crisis happens, the return may be negative but its magnitude is not high. 
The simple return may be unable to record such moment. 
However, we note that, due of the rarity of tailed events, one may use $r_t$ to represent market movement as well.

Because each country's distribution of return differences may have a different average value and variance, we standardize $d_t$ using $z$-transformation; hence, 
\begin{equation}
\label{def:zscore}
	z_{t} = \frac{d_t - \mu_d}{\sigma_d}, 
\end{equation}
where $\mu_d$ is the average of the return difference $d$, $\sigma_d$ is its standard deviation\footnote{$\mu_d$ and $\sigma_d$ are calculated after removing potential outliers. Here, the outliers is defined as any data outside $1.5$ times the interquartile range above the upper quartile and below the lower quartile.} and both statistics are calculated via all observations for consistency. 
Suppose that $D'$ is the set of (country's) return difference after removing outliers. 
Then, the distribution of elements in $D'$ may follow a normal distribution \footnote{see Appendix for verification}. 
The $z$-score in (\ref{def:zscore}) can be then used to separate the tailed events. Since our analysis focuses on the vulnerability of financial network, we create the indicator for volatile capital market. 
The return difference $d_t$ is called \textit{a return jump} when its absolute value of corresponding $z-$score is higher than $2$, otherwise is is called \textit{a noise}.
The direct intuition is that the capital market has approximate $5$ percent probability that it experiences tailed events\footnote{see Appendix for further justification for $c=2$.}. 
\subsection{Probability and Flow}

The structure of the global capital market network is complex.
Standard measures of association, such as Pearson correlation, may be insufficient to reflect the behavioral structure of the financial market network in tailed events.
Financial markets behave asymmetrically; whereas major capital markets (such as those in the United States) can exert significant impact on smaller markets (such as those in Thailand), the opposite is not true.
While the US stock market crash may trigger a global recession, the small stock market collapse is unlikely to do so.
The correlation measure is symmetric and cannot directly capture this relationship.
Additionally, the correlation can best depict the strength of their association when the relationship between two variables is linear.
It then captures the financial market's overall activity rather than just the tails since the tailed events account for approximately 5\% of the whole distribution.
Empirically, the behavior of the capital market's in lower and upper tails may deviate from its normal structure.
The COVID-10 pandemic epitomizes that multiple capital markets plummet simultaneously during the left-tailed events.
The correlation may underestimate the tailed relationship.

One alternative to correlation is the Kendall's tau. 
Its input the two variables's observations and the tau ranges from $-1$ and $1$. 
The positive value indicates the number of concordant observations is higher than discordant ones.
The major advantage of Kendall's tau is that the magnitude of return difference play minimal role, hence it can better capture the tailed events. 
However, its calculation include the entire dataset, not just the tailed distributions. 
Again, it roughly represent the financial market's overall relationship. 
If only tailed events is considered, it may over-exaggerate the degree of association since the left-tailed return difference is strictly less than the right-tailed.
The number of concordant observation will dominate.

We introduce our own measure of return difference's association based on conditional probability of return jumps. 
The (conditional) probability matrix $P$ is an $n$-by-$n$ matrix such that   
\begin{equation} \label{eq:condprob}
 p_{i,j} = \frac{\sum \limits_{t\in T} I_{t,i} I_{t,j}}{\sum \limits_{t\in T} I_{t,i}},
\end{equation}
where $I_{t,i} \in \{0,1\}$ is an indicator that return $i$ jumps at time $t$. 
The value of $p_{i,j}$ is the probability that return $i$ will jump given that return $j$ jumps.
The probability $p_{i,j}$ offers two interpretations; it represent the probabilistic impact of return $j$ on return $i$ and the degree to which return $i$ is dependent on return $j$. 
The higher the probability, the stronger the pairwise dependency structure.
The major advantage of this measure is that $p_{i,j}$ is generally unequal to $p_{j,i}.$
It can then capture the asymmetric behavior of financial markets.
Additionally, it focuses on the behavior at the tailed events by construction. 
The value of $p_{i,j}$ may indicate the strength of tailed relationship rather than the overall behavior. 
Since we attempt to capture the financial system's fragility during the tailed events, we favor the use of conditional jump probability to represent the capital market behavior. 
We note that $p_{i,j}$ capture only probabilistic influence between returns.
It may not indicate the estimated magnitude of return change when one collapses. 
The measure of magnitude of change may be required to fully analyze the financial vulnerabilities.

Our emphasis is on the probabilistic influence between capital returns system, while may entail more than the chance of contemporaneous jumps.
The phrase \textit{flow} from return $i$ to $j$ is used to denote degree of probabilistic impact from return $i$ to $j$.
Our flow formula is inspired by the gravitational forces between two objects, which are proportional to their masses and the distance between them.
The flow from return $i$ to $j$ is contingent upon (1) return $i$'s probabilistic effect on others, (2) return $j$'s ability to accept influence from others, and (3) the \textit{distance} between the two.
The first two factors are nodal and represent the attribute of financial institution or, in our case, country's capital market. 
Despite their nuance, they have different interpretation. 
The first factor is called the \textit{sending mass} and defined for return $i$ as
\begin{equation} \label{def:sendingmass}
	M_i = \sum_{k=1}^n p_{k,i}.
\end{equation} 
The mass $M_i$ captures the overall ability of return $i$ to influence the jump movement in return $j$. 
The greater the value of $M_i$, the more likely that other returns will jump when return $i$ does.
A financial collapse country with high sending mass may correlate with negative jumps in other returns. 
On the other hand, the \textit{receiving mass} of return $i$ is given by 
\begin{equation} \label{def:receivingmass}
	m_i = \sum_{k=1}^n p_{i,k}.
\end{equation}
It denotes the return $i$'s overall capacity to absorb the probabilistic influence of other returns.
The lower value of $m_i$ suggest that return $i$ behaves more independently of other returns (in the tailed distribution).
The sharp return decrease in other countries may be less likely to associate with the change in return of country with low receiving mass.
In this sense, a country with lower $m_i$ can be regarded roughly as influence reducer.  
Although the formula's similarity in (\ref{def:sendingmass}) and (\ref{def:receivingmass}), each country's capital market may have distinct sending, and receiving masses.
Generally, the prominent market has a larger sending mass, whereas the smaller market has a larger receiving mass.

Alternatives to sending and receiving masses are the (weighted sum of) market capitalization and trading volume in the country. 
We note that they generally offer the same interpretation as the sending mass.
A country with high market capitalization and trading volume may have stronger influence on the movement in capital market in other countries. 
Both alternatives to sending mass may be insufficient to represent the market's ability to absorb the probabilistic impact. 
We note that both alternatives should not replace the receiving mass.
Otherwise, the direct inference is that larger capital market will have higher ability to both receive and send out probabilistic influence. 
This may contradict the fact that smaller markets (such as those in Thailand) cannot exert significant probabilistic influence on larger markets (such as those in Europe).

The last factor of flow is the distance $q_{i,j}$ from return $i$ to $j$.
The distance should be asymmetric; $q_{i,j}$ should be unequal to $q_{j,i}$ due to asymmetric relationship in returns.
Additionally, the higher distance between the two markets, the less "connected" their return are. 
The distance $q_{i,j}$ then follows an inverse relationship with the conditional jump probability $p_{i,j}$.
We then use the basic inverse relationship and define 
\begin{equation} \label{def:distance}
	q_{i,j} = \dfrac{1}{p_{i,j}}.
\end{equation} 
The underlying stochastic intuition of $q_{i,j}$ is the average number of times the return $j$ jumps until return $i$ jumps.
The higher $q_{i,j}$ may corresponds to further distance from $i$ to $j$. 
Another important feature of the distance $q_{i,j}$ is that, unlike the first two factors, it reflect the relationship specific to the two market.
For instance, Thailand and Malaysia may have short distance between the two (hence stronger flow) despite their potential low sending and receiving masses. 
The distance $q_{i,j}$ makes it possible to identify the strong relationship in markets with low masses, and weak relationship in markets with high masses.

Our definition of node distance $q_{i,j}$ may be uncommon in the Graph Theory and Network literature.
Typically, a standard hop distance or the shortest path is used to represent distance between nodes. 
This definition may not be applicable in our analysis because our network has high connectivity.
In fact, all markets are connected via the conditional jump probability and the hop distance between any returns pair is then always $1$. 
We discover that some literature create the graph of financial network based on high association measures; two nodes are joined when their respective correlation is higher than certain cutoff values. 
This approach allows different distance between returns in cost of information loss. 
We aim the capture the network structure without removing a relationship from consideration, and exempt ourself from doing so.

Now the formula of flow from return $i$ to $j$ is given by 
\begin{equation} \label{def:flow}
	g_{i,j} = \frac{M_i m_j}{q_{i,j}^2}.
\end{equation}
Again, the formula is similar to the gravitational force between two objects. 
The $g_{i,j}$ captures the overall probability influence of return $i$ on return $j$ and depends on attributes of the two markets.
The flow $g_{i,j}$ is increasing on the sending mass of $i$ and receiving mass of $j$.
When $M_i$ increases, the return $i$ has stronger overall probabilistic influence and the every flow from $i$ increases as a result.
When $m_j$ increases, the return $j$ is able to receive higher amount of flow and every flow into $j$ increases. 
Both nodal attribute shows a positive association with the flow. 
The edge attribute, by contrast, have the opposite direction; the further distance from $i$ to $j$, the smaller respective flow.
We note that the conditional jump probability $p_{i,j}$ also influences the flow $g_{i,j}$ when one inserts (\ref{def:distance}) into (\ref{def:flow}). 
Additionally, Gravity model for international trade also have similar the definition of our ``gravitational" flow in (\ref{def:flow}) because it assumes that the volume of trade between two countries is proportional to their economic mass and a measure of their relative trade distance. 
This model be extended to analyze types of returns to scale by introducing the exponents to each mass and distance. 
We limit our analysis to the constant to scale of flow between returns\footnote{In Appendix, we explores many modifications of (\ref{def:flow}) to study the impact of flow return to scale on our result.}.

\subsection{Shock transmission process}
 
 We adopt our definition of \textit{financial stability} from the global asymptotic stability (GAS).
 Suppose that the financial system is originally at the equilibrium and perturbation occurs. 
 The \textit{stable} system will eventually converges back to the equilibrium, while the \textit{unstable} one will not. 
 In essence, stability is associated with shock resiliency. 
 The unstable capital market network is unable to absorb systemic shocks; the shock magnitude amplifies and follows an explosive path.
 The analysis of financial stability requires investigation of shock transmission process. 
 
 Our shock process construction is heavily inspired by the discussion of \cite{sreejith2017forman} that  
 the Forman Ricci curvature of an edge in a network can be written as the difference between total inflows into and outflows from the edge.
 The concept of flows interaction allows one to capture the behavior of shock transmission in the network.  
 The work from \cite{ricci}, \cite{ricci}, and \cite{2021} show that the Forman Ricci curvature and its alternative measures can be used to create indicators for stability in the financial market.
 Their indicators are simply the average of all curvature values; the indicator rises sharply during instability.
 The flows network and curvatures provide powerful indicators, yet these measure does not have a solid criteria for instability. 
 We aim to develop the instability indicator with additional property that it can systemically distinguish period of vulnerability.
 Hence, we creates the shock transmission process by using flows structure in the capital market network.
 Moreover, because shocks follows explosive nature during instability, our process will scale the magnitude of shock/flow, rather than addition and subtraction.  
 
\begin{figure}
\centering
\begin{subfigure}[b]{0.33\textwidth}
\centering
\includegraphics{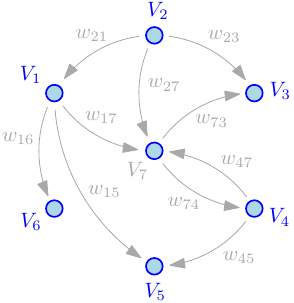}
\caption{Node \\ (Institution)}
\end{subfigure}
\begin{subfigure}[b]{0.33\textwidth}
\centering
\includegraphics{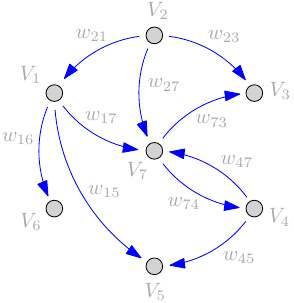}
\caption{Distribution \\ (System complexity)}
\end{subfigure}\begin{subfigure}[b]{0.33\textwidth}
\centering
\includegraphics{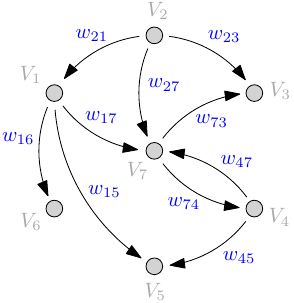}
\caption{Edge \\ (Nodes' relationship)}
\end{subfigure}
\caption{Potential sources of instability.}
\label{fig:shocksources}
\end{figure}
 
Our shock transmission process captures three possible sources of financial network's instability sources, as outlined by figure \ref{fig:shocksources}.
The first source regards the individual characteristics of the node or the financial market. 
When one market experiences extreme shock and crashes, it may spread the the impact and create the cascade of financial collapses.  
On the other hands, capital market in another country may be able to resist the shock and stabilize the system. 
The instability may arise when the majority of nodes amplify shocks. 
The second source is the distribution of shock effects from each market.
It roughly represents the system complexity. 
The more even the shock distribution from one node, the more complex the system.
We note that the complexity may not be indicative of instability. 
The highly connected (hence complex) financial system may be stable because multiple markets absorb the majority of shocks impact. 
The risk-share elements increases the systemic stability.
By contrast, a collapse in one market can lead to multiple breakdowns and system failure. 
This implies that a complex system may be prone to instability from distribution of shocks. 
The last source is the relationship between returns or the (directed) edge in the network.
Like the nodal attribute, edges can create systemic instability when most edges increase the shock impact.
The distinct feature edge characteristic is its asymmetry. 
The directed edge from $i$ to $j$ might increase shock impact while the reverse is not true.
 

The overall shock process is as follows. 
Suppose the financial system is at the equilibrium and systemic shock or small perturbation occurs.
The magnitude of perturbation at each node is scaled by the individual characteristics of the respective market.
Some markets amplify the magnitude, while other reduce.
Afterwards, each market unevenly distributes the shocks to other markets.
The higher the flows from one node to another, the proportion of distributed shock distributed.
These proportions imply the relative dependency between the two market.
After the distribution, the shock from one node to another is further scaled by the (directed) edge-specific factor before arriving at the destination.
Each market receives a new magnitude of total shocks and the process continues until all shocks disappear or become explosive. 

The nodal scaling factor in the first step for return $i$ is given by 
\begin{equation} \label{shock-amp}
	a^N _{i} = \frac{\sum \limits_{k=1}^n g_{k,i} }{ \sum \limits_{k=1}^n g_{i,k} }.
\end{equation}
The multiplier $a^N_i$ is the ratio of total gravitational inflow to outflow for return $i$.
This value is not necessarily one because $G$ is asymmetric.
It is then possible to categorize each return via the scaling factor.  
Return $i$ is called a \textit{shock amplifier} when $a^N_i>1$, 
a \textit{shock absorber} when $a^N_i<1$, and \textit{shock neutral} when $a^N_i=1$.
The value of $a^N_i$ represents how much each return absorbs or amplifies total shocks. 

After node scaling, the total shock on $r_i$ is distributed to return $j\not = i$.
We distribute shocks from $i$ based on the ratio of outflows from $i$. 
The edge $e_{i,j}$ with higher outflows to $j$ receive more shocks from $i$.
Recall that $g_{i,j}$ represent the probabilistic influence of $i$ on $j$. 
The higher distribution ratio shows that shocks will mostly flows in  network with strong relationship.
The ratio of shock distributed from $i$ to $j$ or $u_{i,j}$ is
\begin{equation}\label{predist}
	u_{i,j} = \left( \frac{ g_{i,j}}{\sum \limits_{k=1}^m g_{k,i}} \right).
\end{equation}

The next step is to scale the flow distributed from $i$ to $j$ by the edge scaling factor 
\begin{equation} \label{edgeamp}
	a^E _{i,j} = \dfrac{\text{edge inflow}}{\text{edge outflow}} = 
\left[ \left( \dfrac{w_{i,j}}{\sum_{k} w_{ik}} \right) \sum_{k} g_{ik} \right] / \left[\left( \dfrac{w_{i,j}}{\sum_{k} w_{kj}} \right)  \sum_{k} g_{kj}\right] = \dfrac{\sum_{k} w_{kj}}{\sum_{k} w_{ik}} \cdot \dfrac{\sum_{k} g_{ik}}{\sum_{k} g_{kj}}. 
\end{equation} 
This step is introduced in order to capture the fact that sometimes small influence from one return could trigger paramount impact to another. 
The numerator of $a^E_{i,j}$ represents the estimated outflow from $r_i$ to $r_j$ given $1$ unit of total outflow. 
The denominator, on the other hand, shows estimated the inflow from $r_i$ to $r_j$ given $1$ unit of total inflow. 
The ratio $a^E_{i,j}$ can then approximates the degree of impact adjustment of the flow departing from $r_i$ and arriving at $r_j$.

Let $N$ be an $n$-by-$n$ diagonal matrix such that its $i^{\text{th}}$, $U$ be the outflow distribution matrix, and $E$ be the edge-scaling matrix. Then, it follows that
\begin{equation} \label{shockprocess}
	S_{k} =E_2  \odot (U^T \times N)S_{k-1},
\end{equation}
where $\odot$ and $\times$ denote element-wise and standard matrix multiplication and $k\geq 1$. 
We note that the equation (\ref{shockprocess}) represent our shock process.
The three sources of instability also appear in order as well. 
We note that the stability of the financial network depends on the shock transmission matrix $E_2  \odot (U^T \times N)$.

\subsection{Stability}


Let $n>2$ be an integer and let $\mathbb{A} = \mathbb{R}^+ \cup \{0\}$ be a set of non-negative real numbers.
	Suppose that $E$ and $U$ are $n\times n$ matrices with non-negative entries $e_{i,j}, u_{i,j} \in \mathbb{A}$ for all $i,j$ and that $N$ is a diagonal matrix whose entries are non-negative as well. 
	Define a function $f$ such that 
	\begin{equation} \label{eq:eigen}
	f(\bm{e},\bm{u},\bm{n}) = 	\text{largest eigenvalue of } E \odot (U^T \times N)
	\end{equation}
The shock transmission matrix $E \odot (U^T \times N)$ represents how systemic shocks evolves under current global financial network. 
Its largest eigenvalue denotes the expansion factor of systemic shocks size. 
Note that Perron-Frobenius theorem implies that $f(\bm{e},\bm{u},\bm{n})$ is always positive and its respective eigenvector is also positive as well. 
When $f(\bm{e},\bm{u},\bm{n}) < 1$, the total shock will shrink over time and eventually become negligible. 
The financial network with largest eigenvalue less than $1$ is called \textit{stable} because it can absorb the systemic shock. 
On the other hand, the system with $f(\bm{e},\bm{u},\bm{n}) > 1$ is \textit{unstable}. 

We focus on the impacts of evolution of $E,U$ and $N$ on the dynamics of systemic stability. The analytic implicit function theorem (see, e.g. \cite{KrantzParks}) asserts that the function $f$ is analytic with respect to each component of $\bm{e},\bm{u},\bm{n}$ almost everywhere (as long as the largest eigenvalue is simple).
Hence, we may consider the linear approximation 
\begin{equation}
	f(\bm{e}_2,\bm{u}_2,\bm{n}_2) \approx f(\bm{e}_1,\bm{u}_1,\bm{n}_1) + \sum_{\bm{k} \in \{ \bm{e}, \bm{u}, \bm{n} \}} (\bm{k}_2-\bm{k}_1) \dfrac{\partial f}{\partial \bm{k}} (\bm{e}_1,\bm{u}_1,\bm{n}_1). 
\end{equation}
Under \text{ceteris paribus}, the last term represents the approximate contribution of factor $\bm{k}$ to the total change in largest eigenvalues. 
However, the function $f$ may not be differentiable because our shock transmission matrix is asymmetric positive definite.
The matrices $E$ and $(U^T \times N)$ are not symmetric by construction. 
Thus, we use finite difference to estimate $\dfrac{\partial f}{\partial \bm{k}}$ in our analysis.

\section{Results}
\subsection{Capital market behaviors}

 \begin{figure}
 \centering
 \begin{subfigure}{\linewidth}
 \centering
 	\includegraphics[height =2in]{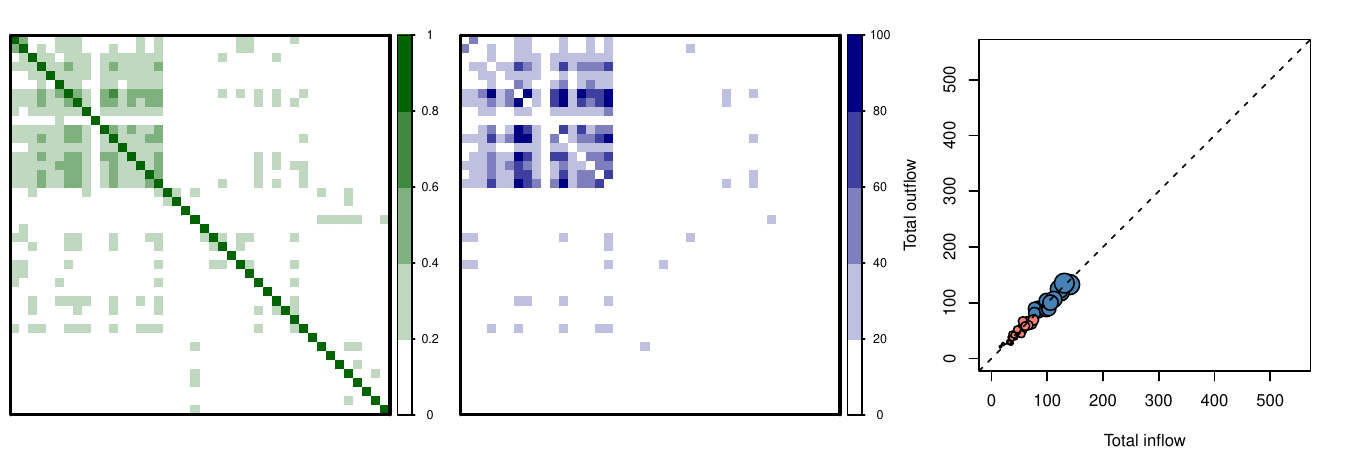}
 	\caption{Overview of capital market structure from January 1995 to July 2021}
 \end{subfigure} 

 \begin{subfigure}{\linewidth}
 \centering
 	\includegraphics[height =2in]{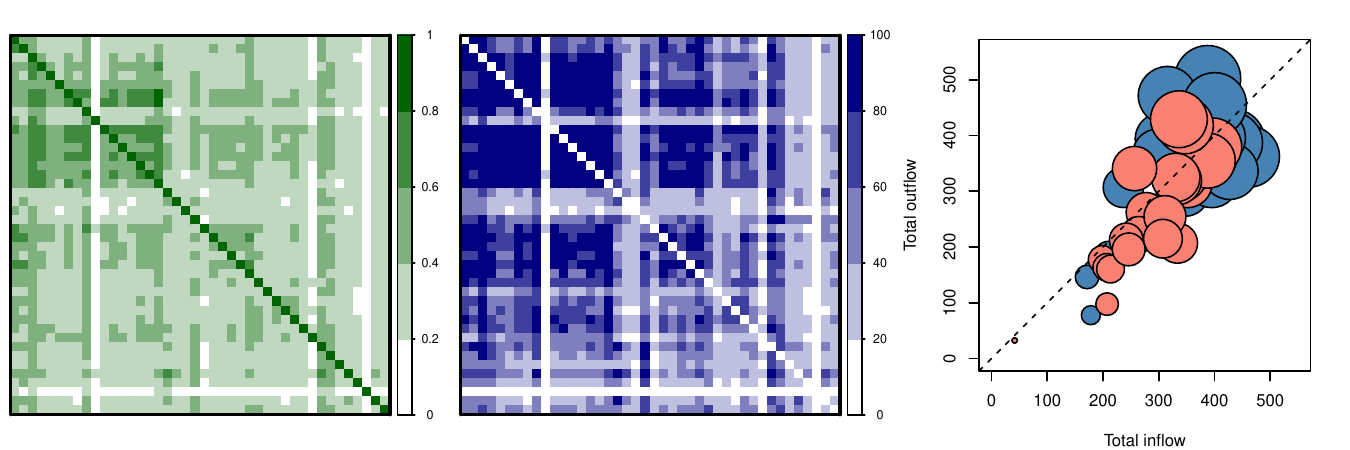}
 	\caption{Overview of capital market structure during global financial crisis}
 \end{subfigure} 
 
  \begin{subfigure}{\linewidth}
 \centering
	\includegraphics[height =2in]{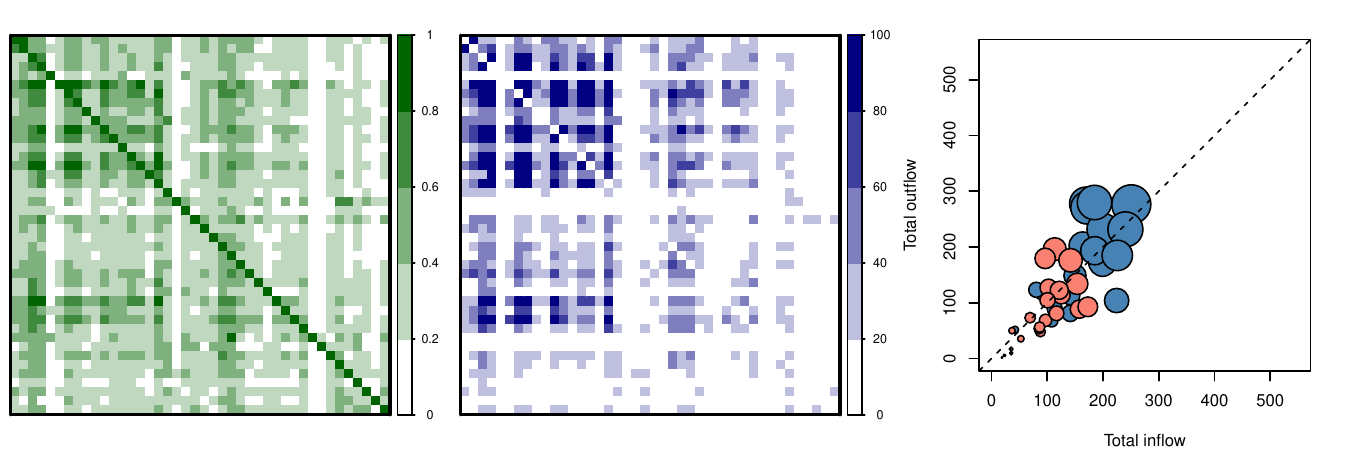}
 	\caption{Overview of capital market structure during COVID-19 pandemic}
 \end{subfigure} 

\caption{The leftmost panel shows the conditional jump probability among 42 capital markets and middle panel shows the flows between markets. The rightmost panel outlines the total inflows and outflows in each market. The size indicates the relative total inflows and outflows. The color red suggests emerging market  and blue suggest developed market.}
\label{fig:cap_overview}
\end{figure}

The behavior of the global financial markets evolves over time.  
Figure \ref{fig:cap_overview} shows the conditional jump probability and the flows structure in the capital markets network across three different period. 
The first period is the entire data set from from January 1995 to July 2021.
It serves to establish the general structure. 
The second and third periods have significant financial and economic events, namely the global financial crisis and the COVID pandemic.
The former covers observations from July 2017 to December 2009, and the latter spans from April 2020 to January 2021.
Both period describe the instability in the capital market. 

The first panel of figure (\ref{fig:cap_overview})  indicates that the global financial network's overall structure is concentrated primarily in advanced economies in North America and Europe.
The top-left quadrant of the conditional probability matrix represents these countries and has a greater value than the remaining quadrants.
This concentration proves that advanced markets in European countries frequently exhibit similarities in their movements.
When one market collapses, the others are likely to follow. 
The dependency structure may be absent in markets outside of the United States; the conditional co-jump probability typically ranges between 0 and 0.2.
As a result, capital markets outside Europe and the United States.
 generally move independently.
 The financial market flows show a consistent structure, but with a greater magnitude. 
The middle heat-map figure demonstrates that the majority of significant flows occur within European countries.
The magnitude of flows in these regions ranges from 40 to 100, and occasionally even higher, while flows elsewhere are less than 20.
The flow structure, hence, emphasizes the network structure that conditional probability captures.
Additionally, the rightmost bubble plot indicates that almost all countries observe roughly equal inflows and outflows, but advanced economies have significantly greater inflows and outflows than emerging markets.
The advanced economy's inflows range between 100 and 200, while the emerging market's inflows range between 0 and 100.
The relative adjusted total flows\footnote{the sum of the average flows in each of six regions} in the market also suggests that after adjusting for regional effects, advanced economy flows remain higher than emerging market flows.

The structure of the capital market network changes abruptly during the devastating financial crises. 
The second and third panels of figure (\ref{fig:cap_overview}) can be used to identify deviations from the overall structure.
The conditional probability of concurrent jumps increases almost everywhere in the probability heat-maps.  
The probability now primarily ranges between 0.2 and 0.6, slightly higher compared from approximately 0 to 0.2 in the general behavior case.
The conditional probability of $0.1$ implies that the average number of times one return jumps until the other follows is approximately $\frac{1}{0.1}=10$, whereas the conditional probability of $0.3$ implies the average number of times is only $\frac{1}{0.3}\approx 3.33$.
A slight increase in probability may correspond to critical change. 
This increase indicates that financial markets have become more interconnected, and a single crisis can spread globally.
Additionally, the flows between markets exhibit a similar structure.
Each country's inflows and outflows increase significantly.
Additionally, the majority of countries receive more inflows from, send more outflows to other countries.
Therefore, the flow network is more connected than in the general case.
We remark that the flows during the COVID pandemic may seem less connected.
This is because the time at which each country experiences severe outbreak is different. 
Although the capital market return during the outbreak is extremely volatile, the volatility decreases when the epidemic is controlled. 
This means the conditional jump movement may not be as high despite the fact that the pandemic affects every market. 
Additionally, each country's total inflows and outflows differ significantly from the corresponding general behavior case.
Throughout the global financial crisis, the majority of inflows and outflows have been between 200 and 550.
This increase is nearly three times the size of the initial range.
Surprisingly, the majority of advanced economies experience greater inflows than outflows, whereas emerging markets experience greater outflows than inflows.
One possible conclusion is that the most advanced economies amplify the size of flows, thereby magnifying the severity of shocks or perturbations.
On the other hand, emerging markets mitigate the systemic shocks magnitude.
The increase in the inflows and outflows during the COVID pandemic follows the similar intuition, but with less magnitude.

\subsection{Financial fragility indicator}

The evolution of the financial network stability can be analyzed through the largest eigenvalue of the shock transmission process in (\ref{shockprocess}) across different time span $T$.
Our time span set will be a $120$-day rolling window starting from July 1995 to July 2021. 
There is no general consensus of the appropriate window size, but we believe the market behavior in approximately 6-month is suitable for analysis.

\begin{figure}[h]
\centering
\includegraphics[width=\linewidth]{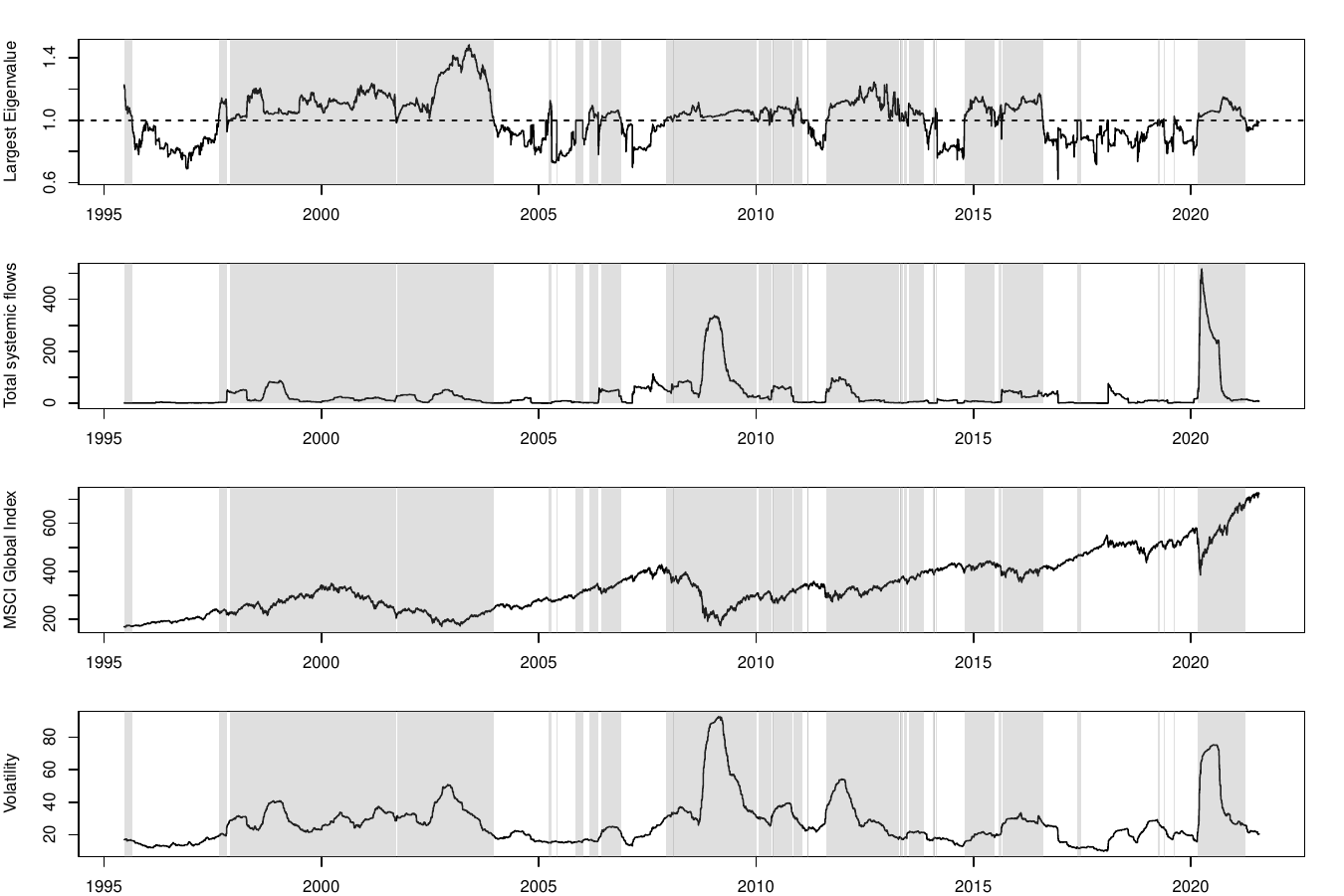}
\caption{The first panel shows the largest eigenvalue of the shock transmission process can be indicative of the systemic stability. Its value at the stable is below $1$. The shaded area indicates the period of instability. The remaining panels shows related statistics. The total systemic flows is the total amount of inflows and outflows in the system and can be used to proxy the magnitude of shock impact. The MSCI global index and return serve to verify the robustness of the result. }
\label{fig:fragility_index}
\end{figure}


The first panel of figure \ref{fig:fragility_index} demonstrates that there are five significant period of instability, each seemingly connect to major global financial events.

The first instability period occurs between August 22, 1997 and December 16, 2003. 
It is the longest period of instability and is characterized by a series of financial events.
In the summer of 1997, the Asian financial crisis began, spreading panic throughout the capital market network.
The crisis ended in 1999 and multiples of Southeast Asian countries suffers from severe economic depression and stock market collapses. 
In 2000, the United States's technology bubble burst and the so-called "Dot.com" recession began in 2000.
Then, the September 11 attacks in 2001 also sparks the panic in the global financial network.
The US stock market crashed again in 2002 because of the Enron corporate fraud and loss in investor confidence. 
Additionally, geopolitical tensions between the US and the Middle East, as well as the 2002-2004 SARS outbreak, may contribute to the instability.
Despite the fact that this period of instability is extremely long and the largest eigenvalue is relatively large, the total systemic flows over the entire time span are relatively small.
This implies that the underlying capital market structure is not highly connected, and thus many markets may behave autonomously.
One market's impact on another may be negligible.
Additionally, the MSCI global index conforms to intuition.
The index briefly fell during the Asian financial crisis of 1998, but then resumed its upward trend.
This is consistent with the bursting of the technological bubble.
After the late 2000 bubble burst, the index began a downward trend, paralleling the critical financial crash.
Additionally, the MSCI global index exhibits relatively high volatility during the period, compared to those prior to 1997.

The second instability period correspond to the global financial crisis in 2007-2011.
It began in December 4, 2007, consistent with housing bubble burst in the US equity market and ended at March 8, 2013.
Surprisingly, the largest eigenvalues in this period only stay slightly above $1$.
Such values indicate stability, but not the total impact within the financial system. 
By contrast, the total systemic flows soars and reached the highest point during late 2008 to early 2009.
This is just after Lehman Brothers bankruptcy and its end after the announcement and implementation of quantitative easing policy in various countries such as the US and UK.
We note that, when the total systemic flows reached at peak, the MSCI global index plummets and the the respective return becomes extremely volatile.
This peak may mark the lowest point in the crisis. 
Additionally, the periods has significantly higher amount of flows than the previous ones. 
Hence, the financial network has become more connected and one collapse in one market will likely correspond to multiple collapses in other markets.
Following up as the aftermath of the global financial crisis is the European debt crisis in 2011-2014. 
This third instability period occurred between from August 18, 2011 to February 28, 2014. 
The major event here is Greece's government debt crisis began and the panic spreads again.
The total systemic flows here is not as high as in 2008, but bear considerable magnitude. 

 The fourth period of instability is primarily focused on the legal process surrounding Brexit.
 The period covered is October 20, 2014 to August 8, 2016.
Although the largest eigenvalues are typically greater than $1$, their fluctuation provides stability for a brief period of time.
Total systemic flows are surprisingly low in comparison to other periods of instability.
Perhaps the financial system will be less interconnected following Brexit.
This suggests that the instability associated with Brexit may be milder than that associated with the global financial crisis, as the systemic impact is likely to be much smaller.
Hence, it is possible to observe movements in capital markets outside the United Kingdom, but economic recessions.
Another unique feature of this fragility period is that the MSCI index moves slightly and hovers around 400 in the whole duration.
The respective MSCI return also suggests a low level of volatility.
Perhaps, the global capital markets are unsure what to make of Brexit's impact and progress. 
 
The most recent period of instability is associated with the COVID-19 pandemic.
It is the shortest and, most likely, the most severe period.
The outbreak began on April 21, 2020, when most countries implemented a lockdown policy, and it has since spread throughout the world.
Total systemic flows increase abruptly to a historical high before gradually decreasing.
The sudden increase suggests that the financial markets became suddenly highly connected to one another.
It is unsurprising that capital returns have collapsed almost everywhere, and that most economies have entered a recession.
At the start of the instability, the MSCI index also shows level shifts.
This abrupt change indicates that the pandemic is unexpected, and the corresponding market return is also extremely volatile.
However, we are struck by the extreme speed with which the market is adjusting.
Around August 2021, systemic flows return to near zero.
The MSCI also returns to pre-crisis level, and volatility falls slightly.

\subsection{Contribution to Stability Analysis}

We continue our investigation by calculating the effects of the potential instability sources depicted in Figure 1 on the change in the largest eigenvalue.
The contribution to stability is split into two categories.
The first scenario occurs when the system approaches stability and the largest eigenvalue decreases.
The value does not have to be greater than 1.
The second case is when the financial network becomes unstable and the largest eigenvalue rises.
It may seems counter-intuitive that the positive contribution to change in largest eigenvalue implies higher instability.
This is because increase in largest eigenvalue implies instability. 

\begin{figure}[h]
\begin{subfigure}{\linewidth}
	\includegraphics[width=\linewidth]{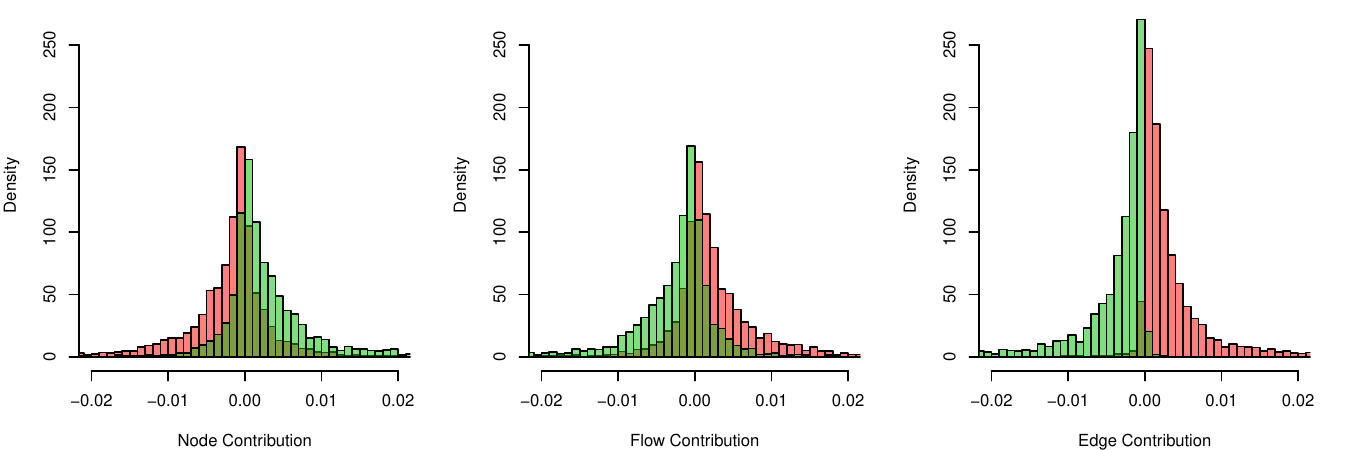}
	\caption{Histogram of contributions to change in the largest eigenvalue. }
\end{subfigure}

\begin{subfigure}{\linewidth}
	\includegraphics[width=\linewidth]{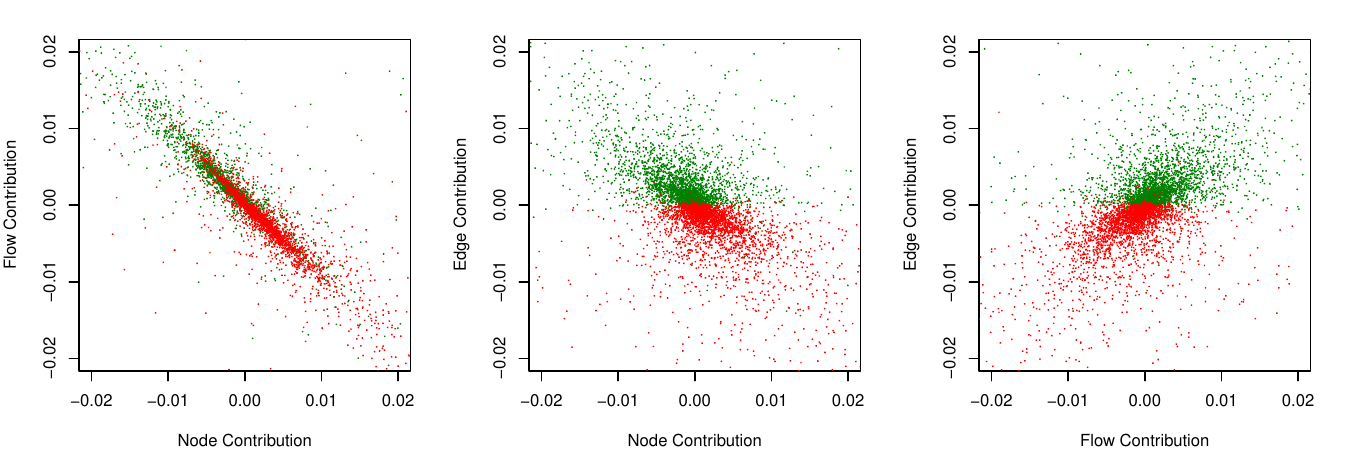}
	\caption{Associations between contributions to change in the largest eigenvalue.   }
\end{subfigure}

\caption{The contributions to stability and instability of nodes, flows, and edges. The green component indicates the contributions when the financial network moves towards stability and the largest eigenvalue decreases, while the red suggests the opposite.}
\label{fig:fragility_eigen_contr}
\end{figure}

Figure  \ref{fig:fragility_eigen_contr} depicts the distribution of each shock process component's contribution to change in stability. 
The contribution of node or individual markets to stability appears to be inversely proportional to the change in largest eigenvalue. 
When the system is approaching stability, node contributions are generally positive, with most ranging from 0 to 1 percent.
Changes in individual capital markets may, in turn, cause systemic instability.
However, when the financial network becomes unstable, the node contribution is negative.
This means that the financial stability is created by the combination of the individual characteristics of each capital market.

The effects of market individual characteristics appeared to be offset by flow distribution effects.
Unlike that of nodes, the contribution to the stability of the flow distribution is proportional to the change in the largest eigenvalue.
Hence, the distribution of shocks influence stability, and the magnitude is close ti that from individual country.  
Surprisingly, the contributions of node and flow have strong negative associations.
When one rises, the other falls.
This may explain the significance of systemic adaptation to (in)stability.

Assume that the capital markets become unstable as a result of a financial crisis in one country.
Such markets act as shock absorbers, increasing the respective node-scaling factor and connecting the financial system.
An increase in the node-scaling factor may promote stability on its own.
This is because it implies that the corresponding market will have a greater impact.
When the shock distribution does not adjust and the financial system is not highly connected, the financial crisis may spread to a small number of markets.
The system will be able to absorb shocks, and its stability will improve.
However, as the network becomes more connected during the crisis, the adjustment in flow distribution implies greater instability.
The impact of one country quickly spreads around the world.
It is unclear what the combined effects of flows and node distributions will be.

The contribution to edge stability or market relationship is directly proportional to the change in largest eigenvalue.
When the system is approaching stability, the contribution of the edges is negative and ranges (again) from 0 to 1 percent. 
There is little probability of contribution in the ``wrong" sign.
This implies that market relationships may be a major driver of systemic stability.
However, it also exhibit decent association nodes and flows contribution, possibly because of network complexity during instability.
The shock impact from one node spreads more evenly to other markets as the financial system becomes more connected.
Shock outflows to multiples markets receives more inflows and the edge-scaling factors play more role determining total systemic shock size.
As a result, when the flows contribution is positive, the edge contribution is likely to be positive as well.

\subsection{Connection to volatility}

Apart from the analysis of financial stability from the largest eigenvalue, figure \ref{fig:fragility_index} also shows that the behavior of the volatility of the global MSCI index return and the total systemic flows are almost identical. 
Therefore, the total systemic flows can be used as the proxy or alternative for systemic volatility.
The major advantage of this measures is the tractability; it is possible to identify flows contributing to the volatility.

\begin{figure}[h]
\begin{subfigure}{\linewidth}
 \centering
    \includegraphics[width = \textwidth ]{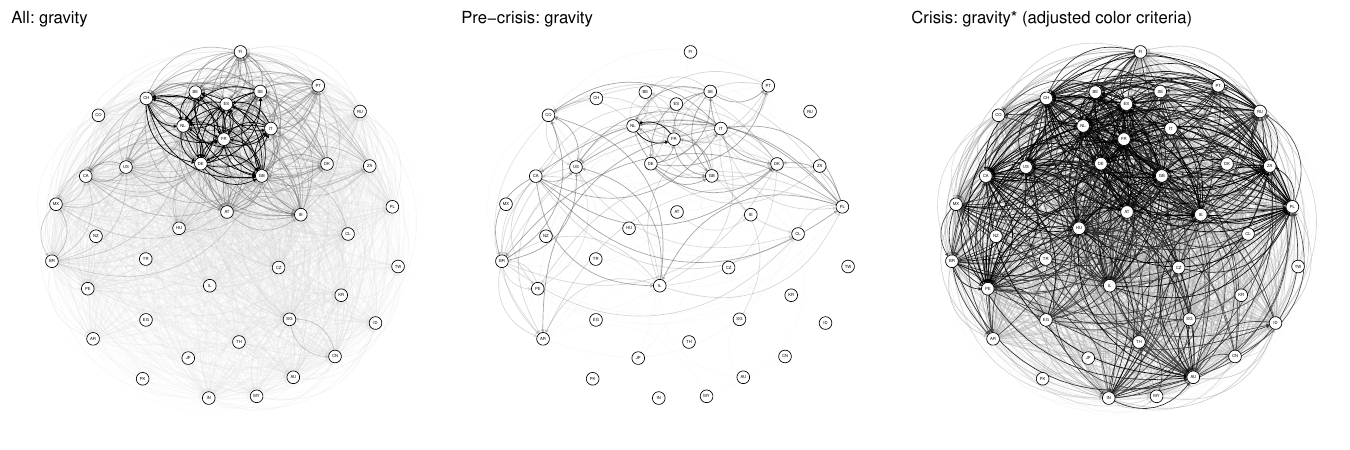}
    \caption{Flows in the capital market before and during the COVID-19}
\end{subfigure}

\begin{subfigure}{\linewidth}
    \centering
    \includegraphics[width = \textwidth ]{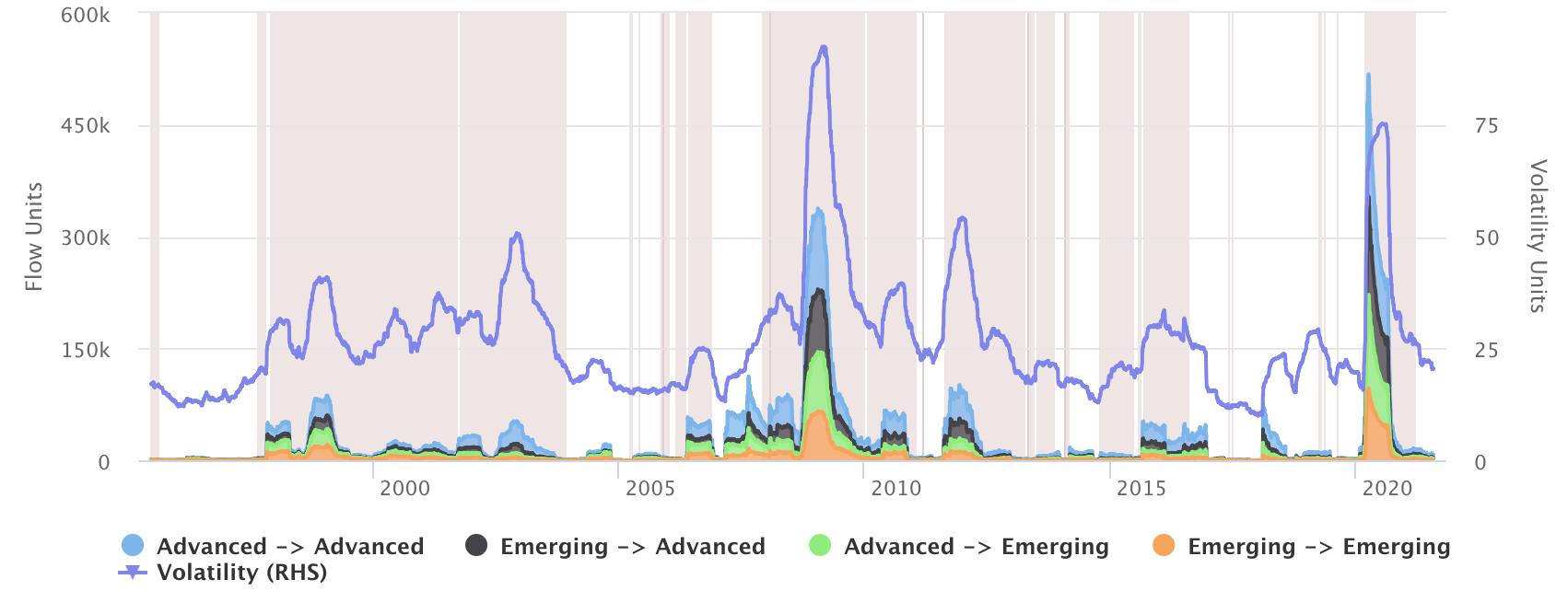}
    \caption{Financial market volatility and decomposition of total flows by market types}
\end{subfigure}

\caption{Market volatility and the total flows in financial network shows almost identical movement. Shaded area show instability periods determined in figure \ref{fig:fragility_index}}
\label{fig:scenario_gravity}
\end{figure}

The magnitudes of flows in the capital market network are depicted in panel (a) in figure \ref{fig:scenario_gravity}. 
The first network indicates that financial market volatility is generally caused by flows within the advanced economies of Europe and North America. 
Flows within and around such a region are stronger than in other areas. 
As a result, when one country experiences critical financial events and its return soars, the other countries may follow. 
Concurrent return jumps are more likely to spread within the region and less likely to spread to emerging Asian markets. 
Volatility movements in general may thus be influenced by developments in advanced capital markets in the United States and Europe. 
On the other hand, the magnitude of flows in emerging markets, particularly in Asia, is significantly lower. The region's inflows and outflows from Singaporean and Chinese markets are comparatively stronger, albeit of minor magnitude. 
This means that the volatility of capital markets in such regions is influenced by these two countries, which may be influenced by the European market. 
They may act as a buffer between the Asian and European markets.

The flows within the financial market network prior to the COVID-19 pandemic demonstrate a shaky relationship between markets. 
The calculation includes market returns in the 120 working days preceding January 1, 2020. 
There were no pivotal financial events during these periods, and the magnitudes of flows were smaller. 
The flows are mostly within European and North American countries, indicating the sources of volatility. 
We observe that, due to the size of the flows, the conditional probability of concurrent jumps provided by one return jump is significantly lower. 
With a few exceptions, the movement of capital markets in this period may be generally independent on most countries.

During the COVID-19 crisis, the financial market's behavior changes abruptly. Figure \ref{fig:scenario_gravity}'s rightmost network depicts the flows in 120 working days prior to April 1, 2021, just three months after the behavior described in the middle panel. 
The magnitude of flows everywhere surges ominously. 
Rather than a single country or region, every market appears to be the source of systemic volatility. 
The market is extremely fragile, and the underlying dependency structure is extremely strong. 
When the magnitudes and spreads of the flows are larger globally, the conditional probability that two returns will jump simultaneously when one greatly increases or decreases is significantly higher than before. 
The fall in one country's return may set off a chain reaction of subsequent return collapses in other countries. 
Furthermore, the magnitude of the impact is greater. 

Panel (b) of figure \ref{fig:scenario_gravity} shows the breakdown of total network flows by market type (advanced and emerging). 
The total flows and market volatility shows extremely similar movement; hence, the former can be used as the proxy for the latter. 
It is then possible to determined the main source(s) of market volatility. 
Roughly speaking, the flows form advanced markets contribute approximately about 50 to 75 percent of the market volatility. 
However, such total contributions decreased a bit during instability periods; the contribution of flows from emerging market is almost the same as that from advanced economy.
These two observations may suggest that the main volatility source is from the advanced market.
The movement or disruption in an advanced market can generally triggers volatile nature of the global financial network. 
During instability, almost every market may co-move together and contribution of flows from emerging market increases as a result. 
This also partly reflects high interconnectedness of the financial network at such periods as well.
\section{Conclusion}

The analysis on the flows network shows promising result and future improvements are possible.

The flows between capital market's returns is able to capture the market behavior. 
The overall total flow magnitude is stronger when the financial market enters the unstable state.
Additionally, the pairwise flow shows higher magnitude when two returns share stringer movement. 
The evolution of flows between market can therefore be used to analyze and compare network's behavior and to signal incoming critical transition into instability or stability. 
Although it cap capture the probabilistic impact of one return on another, the underlying calculation does not involves the estimated impact magnitude.  
This is current beyond the scope of our analysis at the moment. 
However, should the expected increased or decreased in returns is considered in the flows calculation, the resulting network can then be used in stress test.
The impact of each financial crisis can then be quantified and analyzed. 

The financial stability indicator constructed from largest eigenvalue of the shock transmission network also correspond to empirical observation. 
The suggested instability period generally corresponds to the ones where the global MSCI index falls or the financial market experience significant events like Brexit. 
The strength of this indicator lies in its stability criterion; most alternative instability indicator only shows abrupt movement during instability while not having solid criteria to separate stability. 
It is also possible to calculate the contribution to change in the stability indicator and determined the main factor(s) that creates instability.  
This approach can be useful in an analysis to determine the fragile locations in the financial network. 
We note that the movement of our indicator is quite volatile due to the stochastic nature of capital market returns in the short run.
When added more days into the rolling calculation, the indicator becomes much smoother, but loses some of its predictive power.
The appropriate length of the rolling window can be further explored. 
  
Surprisingly, the movement of the market volatility and total network flows share almost identical movement despite different calculation methods.
The market volatility is the annualized standard deviation of the global MSCI index, while our flows are determined from conditional probability. 
We believe that the total network flows can be alternative measures of market volatility.
The advantage of this measure is the tractability; one can determine the major contributors of volatility from our approach.
We found that the total flows from advanced market mainly contribute to market volatility.
More analysis could be performed to determined the possible sources of volatility.
It would greatly benefit portfolio optimization and analysis on pre-emptive policy. 
We note that more analysis could be performed to extract information on the nature of volatility.

Our network flows approach can disentangle the interconnectedness within and stochastic nature of the global financial market. 
Because the financial system has always been evolving, today's network structure today may not be representative of tomorrow's.
Our approach may be able to determine the today's sources of volatility and instability.
If provided some views on the future market movements, it is possible to construct another network to explore its nature and stability. 
Such modifications can help today's analysis on estimated the impact of possible future change and help investor and policymakers to potentially react and prepare.

\clearpage
   
\printbibliography

\clearpage

\section{Appendix}

\subsection{Verification of cutoff value $c=2$}

After removing outliers, most capital market returns, except Egypt's, follow a normal distribution.

\begin{figure}[h]
\centering
	\includegraphics[width = 5in]{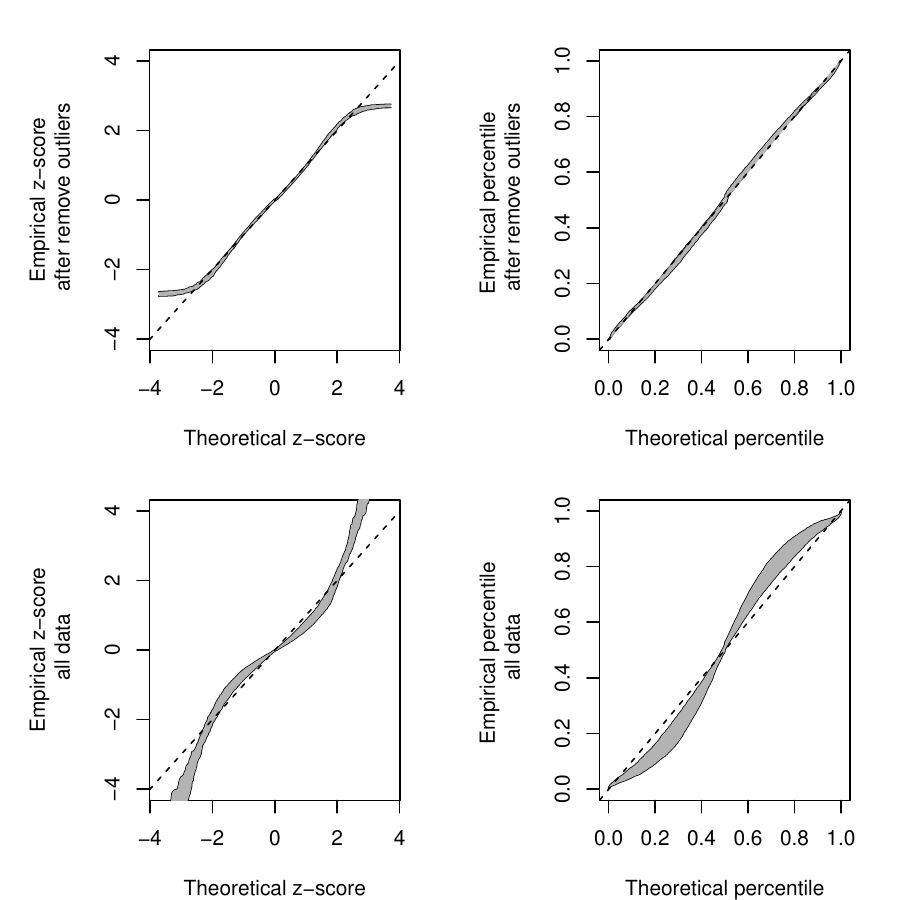}
	\caption{After removing potential outliers, most returns follows a normal distribution. However, the returns behave differently when the absolute value of underlying $z$-score is higher than $2$. The shaded area represents the range of the P-P and Q-Q plot across different returns. }
	\label{fig:robust_cutoff}
\end{figure}

The top panel of figure \ref{fig:robust_cutoff} shows that the normal distribution fits the return difference data after removing outliers almost perfectly. 
The theoretical and empirical percentiles are almost equal. 
However, the theoretical percentile may overestimate and underestimate the empirical behavior when the respective theoretical $z$-score is higher than $2$ and below $-2,$ respectively.
Thus, the cutoff value of $2$ seems appropriate to identify the return jumps. 
More over, the underlying interpretation of $z$-score outside the interval $[-2,2]$ is the percentile approximately below $2.5$ or above $97.5$.
The estimated percentage of outliers/jumps in the return difference distribution is then about five percent. 

The typical normal distribution fails to represent the return difference's behavior. 
When the theoretical percentile is below (above) $50$ percent, the normal distribution generally overestimates (underestimates) the empirical observations.
The estimate at the extreme value with $z$-scour outside $[-2,2]$ also shows the opposite direction.

\subsection{Robustness check}
We calculate different financial stability indicators and compare across our baseline provided in previous section. 
There are five aspects for robustness check. 

\begin{figure}
\begin{subfigure}{\linewidth}
	\includegraphics[width = \linewidth]{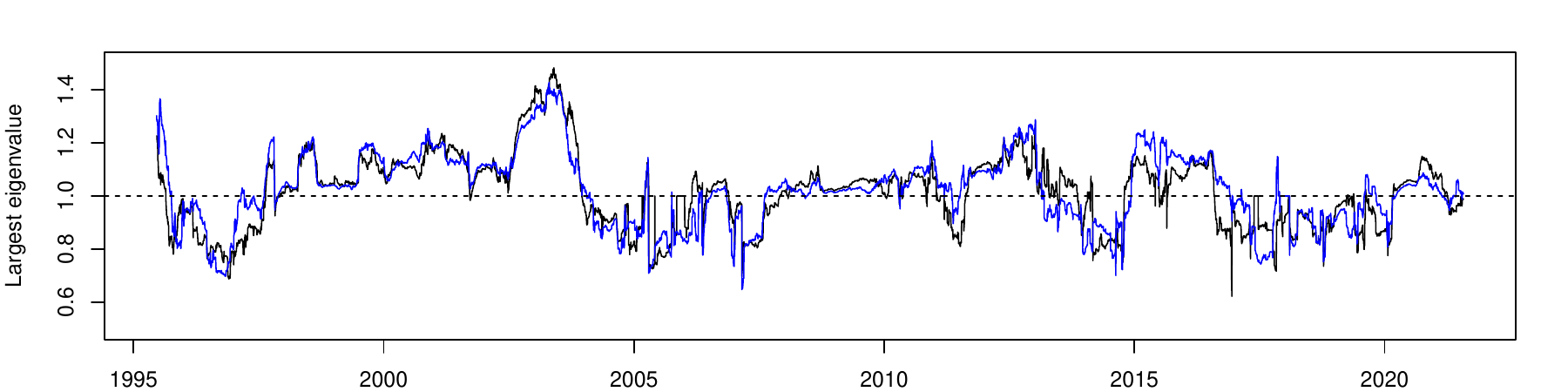}
	\caption{ Replacing return difference with simple return in (\ref{def:zscore})
	}
\end{subfigure}

\begin{subfigure}{\linewidth}
	\includegraphics[width = \linewidth]{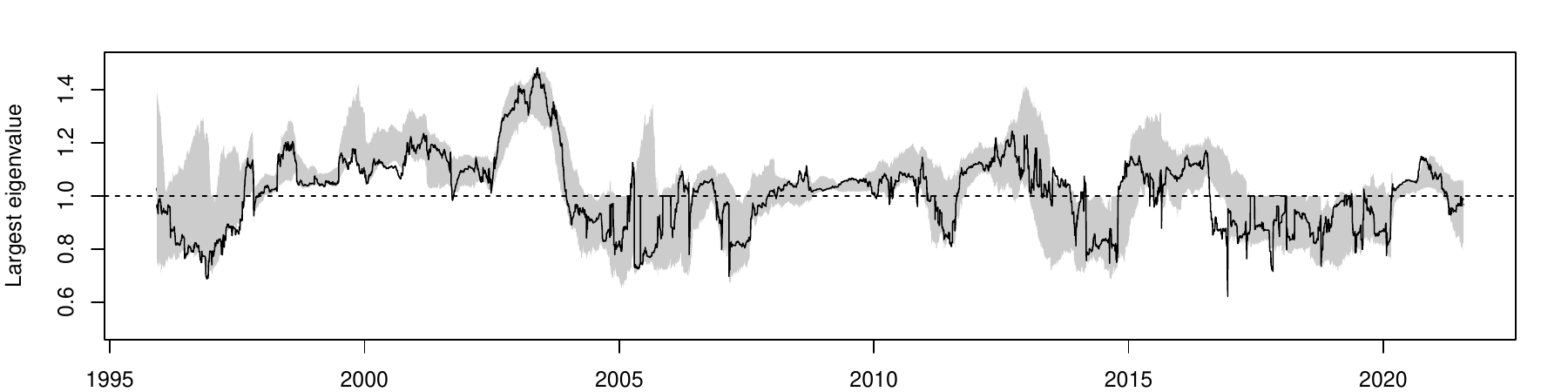}
	\caption{Changing different time-shift $\triangle t = 1,2,3,...,120$}
\end{subfigure}

\begin{subfigure}{\linewidth}
	\includegraphics[width = \linewidth]{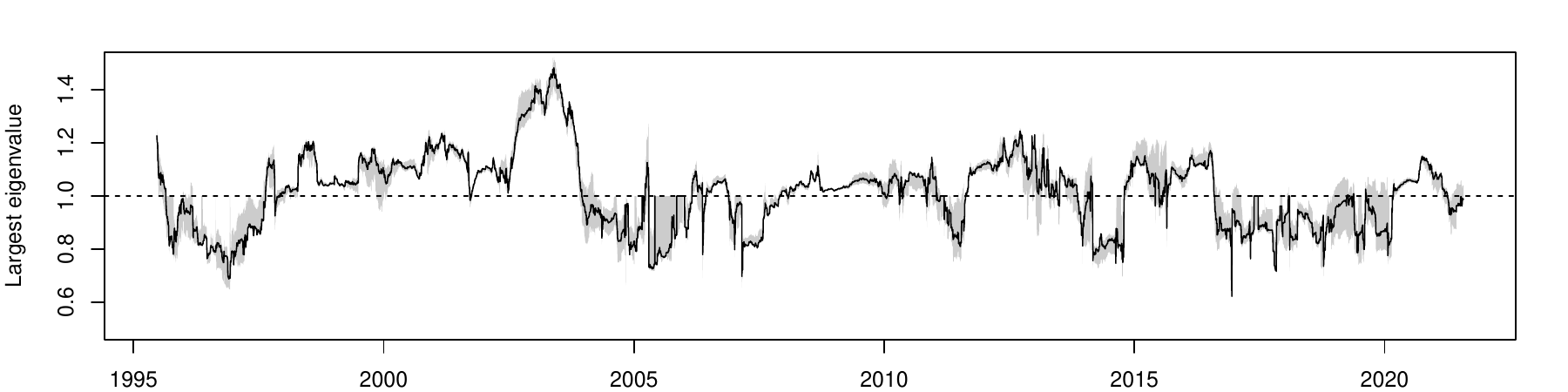}
	\caption{Changing different cutoff $c$ from $1.8$ to $2.2$ with $0.005$ increment}
\end{subfigure}

\begin{subfigure}{\linewidth}
	\includegraphics[width = \linewidth]{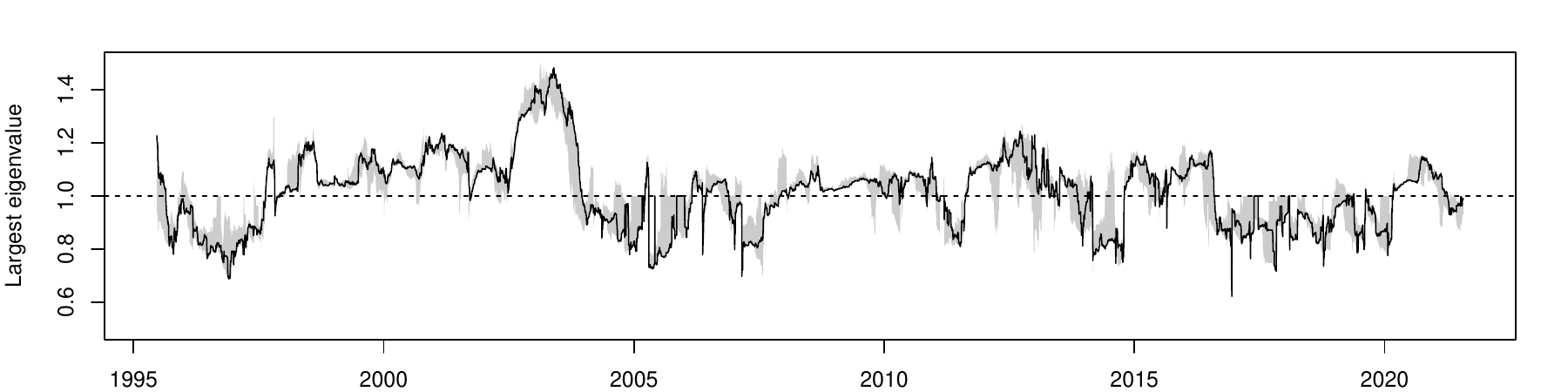}
	\caption{Changing rolling window $60$-,$61$-,$62$-,...,$120$-day. }
	\end{subfigure}
\caption{Robustness check from different aspects. The black line indicate the baseline fragility indicator. The blue line indicates the alternative created by using  simple return in (\ref{def:zscore}). The shaded area shows the $5^{\text{th}}-95^{\text{th}}$ range of the all alternative indicators in each aspect.}	
\label{robust:all}
\end{figure}


The first aspect is about variables that represent the abrupt change in the behavior of the capital market. 
Our analysis uses the return difference to determine the periods with high volatility in equation (\ref{def:zscore}). 
The high return difference implies that current return increases or decreases greatly from the previous day's, hence indicating significant financial events.
It is possible to use the actual return to identify periods with abrupt change as well.
This is because most returns have an average of $0$.
Therefore, the distribution of return difference, and the value of return will have an average of $0$, but their corresponding variance will differ.
Replacing return difference $d_t$ with return $r_t$ in  (\ref{def:zscore}) will provide different interpretation of period of volatility.
The significant financial events happens when the underlying capital market index greatly drops or rises. 
The top panel of figure \ref{robust:all} shows two financial stability indicators constructed from two approaches.
The black indicator black is our baseline version, and the blue one shows the alternative.
The alternative indicators is constructed by replacing return difference $d_t$ with return $r_t$ in  (\ref{def:zscore}).
Calculations in (\ref{eq:condprob})- (\ref{eq:eigen}) remain the same.
It follows that both indicators share similar trend and produce almost identical period of instability.

The second aspect concerns the use of time-shift factor $\triangle t$ to calculate capital market return.
Our analysis set  $\triangle t=1$ to analyze daily return.
Variable $\triangle t$ can be any positive integer theoretically, and the standard value can be $5,20,$ and $60$ for weekly, monthly, and quarterly return respectively. 
An increase in the time-shift value should emphasize changes in the index; 
however, the return calculated from higher value of time-shift may fail to include or reflect financial events within the corresponding period.
Foe example, the yearly return (time-shift approximately $240-255$) may fail to capture the financial events in the first two months. 
We believe that appropriate value of $\triangle t$ is arguably between $1$ and $120$ inclusive.
The second panel shows the $90$-percent interval of the fragility indexes created by varying the time-shift factor $\triangle t$ in (\ref{def:return})
All indexes show similar general trend, and their implied instability periods are similar.
We note that when the interval has quite high range when the baseline indicator suggests stability.
This may suggests that the stability period may have hidden instability.
On the other hand, the interval is narrow and hover above $1$ when the baseline indicator suggests instability.

The cutoff value in equation (\ref{def:zscore}) is also vital to our stability index.
We note that our result is sensitive to the cutoff value $c$.
This is because $c$ is used to distinguish periods with  significant financial movements.
The value too low may include too much noise, while the value too high may not include sufficient data.
Because the returns difference after removing outliers follows normal distribution, we believe the cutoff value from $1.8$ to $2.$ is appropriate.
The standard normal distribution with when $z-$score equals to $1.8$ and $2.2$ have $7.2$ and $2.8$ percent of potential outliers, respectively.
The third panel of figure \ref{robust:all} shows the middle $90-$percentile range of the stability indicators calculated by altering the cutoff $c$ in  (\ref{def:zscore}).
The range is extremely narrow, and has the same trend as the baseline indicator. 

Another important aspect is the rolling window size of the sub-period set $T$.
It serves to capture the behavior of capital market networks within the specific time period.
The lower the window size, the more recent the behavior captured.
However, when the window size is too small, the number of observations will be insufficient to properly reflect the financial market structure.
The resulting financial stability indicator may be volatile and fail to provide useful insights.
On the other hand, when the window size is too large, the behavior of the financial market may be dominated by the observations far in the past.
This means it may not reflect the current behavior of the market.
We believe that the window size from $60$ to $120$ days are appropriate for analysis. 
It should be able to capture the financial market behavior within $3$ and $6$ months, respectively.
The last panel of figure \ref{robust:all} also shows the that the financial fragility indicators produced from various rolling window size provides similar behaviors and are almost identical to each other.
\end{document}